\documentclass{svmult}%
\usepackage{makeidx}
\usepackage{graphicx}
\usepackage{multicol}
\usepackage[bottom]{footmisc}
\usepackage{natbib}
\usepackage{url}
\usepackage{amsmath}
\usepackage{amsfonts}
\usepackage{amssymb}%
\renewcommand{\PackageWarningNoLine}[2]{}

\begin{document}
%
\title*{Coarse spaces over the ages}
\author{Jan Mandel\inst{1} \thanks{
Supported by National Science Foundation under grant DMS-0713876.} \and
Bed\v{r}ich Soused\'{\i}k\inst{1,2} \thanks{
Partially supported by National Science Foundation under grant DMS-0713876
and by the Grant Agency of the Czech Republic under grant 106/08/0403.} }
\institute{Department of Mathematical and Statistical Sciences, \\
University of Colorado Denver, Campus Box 170, Denver, CO 80217, USA
\and Institute of Thermomechanics, Academy of Sciences of the Czech Republic,
Dolej\v{s}kova 1402/5, 182~00 Prague~8, Czech Republic.
\texttt{jan.mandel@ucdenver.edu, bedrich.sousedik@ucdenver.edu}}
\maketitle

\section{Introduction}

\label{mandel_mini_3_sec:introduction}

The objective of this paper is to explain the principles of the design of a
coarse space in a simplified way and by pictures. The focus is on ideas rather
than on a more historically complete presentation. That can be found, e.g., in
\cite{mandel_mini_3_Widlund-2008-DCS}. Also, space limitation does not allow
us to even the mention many important methods and papers that should be
rightfully included.

The coarse space facilitates a global exchange of information in multigrid and
domain decomposition methods for elliptic problems. This exchange is
necessary, because the solution is non-local: its value at any point depends
on the right-hand-side at any other point.
Both multigrid and domain decomposition combine a global correction in coarse
space with local corrections, called smoothing in multigrid and subdomain
solves in domain decomposition. In multigrid the coarse space is large
(typically, the mesh ratio is 2 or 3 at most) and the local solvers are not
very powerful (usually, relaxation). In domain decomposition, the coarse space
is small (just one or a few degrees of freedom per subdomain), and the local
solvers are powerful (direct solvers on subdomain). But the mathematics is
more or less the same.

\section{Local nullspace and bounded energy conditions}

Consider the variational problem
\begin{equation}
u\in V:\qquad a(u,v)=f(v)\qquad\forall v\in V,
\label{mandel_mini_3_eq:variational}%
\end{equation}
where $a$ is symmetric positive definite and $V$ is a finite dimensional
space. Most, if not all, multigrid, domain decomposition, and substructuring
methods for (\ref{mandel_mini_3_eq:variational}) can be cast as variants of
the additive Schwarz method (ASM), which is the preconditioning by the
approximate solver
\begin{equation}
M:r\mapsto\sum_{i=0}^{N}u_{i} \label{mandel_mini_3_eq:additive-prec}%
\end{equation}
where $u_{i}$ are solutions of the subproblems
\begin{equation}
u_{i}\in V_{i}:\qquad a(u_{i},v_{i})=r(v_{i})\qquad\forall v_{i}\in V_{i}
\label{mandel_mini_3_eq:subproblem}%
\end{equation}
where%
\begin{equation}
V=V_{0}+V_{1}+\cdots+V_{N} \label{mandel_mini_3_eq:space-dec}%
\end{equation}
The resulting condition number of the preconditioned problem is then bounded
by $nC_{0}$, where $n\leq N+1$ is the maximal number of the subspaces
$V_{0},V_{1},\ldots,V_{N}$ that have nontrivial intersections, and $C_{0}$ is
the constant from the bounded energy decomposition property%
\begin{equation}
\forall v\in V\exists v_{i}\in V_{i}:v=\sum_{i=0}^{N}v_{i},\quad\sum_{i=0}%
^{N}a\left(  v_{i},v_{i}\right)  \leq C_{0}a\left(  v,v\right)  ,
\label{mandel_mini_3_eq:energy-bound}%
\end{equation}
cf.,
\cite{mandel_mini_3_Widlund-1988-ISM,mandel_mini_3_Bjorstad-1991-SSO,mandel_mini_3_Dryja-1995-SMN}%
.

Variants of ASM include the multiplicative use of the subspace correction in
\cite{mandel_mini_3_Mandel-1993-BDD,mandel_mini_3_Mandel-1994-HDD}, and the
use of other forms in place of $a$ in subproblems
(\ref{mandel_mini_3_eq:subproblem}), cf.,
\cite{mandel_mini_3_Dryja-1995-SMN,mandel_mini_3_Toselli-2005-DDM}.

Now consider $V$ to be a space of functions on a domain $\Omega$. The
subspaces $V_{i}$ range from the span of one basis vector in multigrid (for
the simplest case, Jacobi iteration) to spaces of functions on large
overlapping subdomains $\Omega_{i}$. When the domain $\Omega$ is the union of
non-overlapping subdomains $\Omega_{j}$, $j=1,\ldots,M$, the spaces $V_{i}$
are defined as certain subspaces of the space $W=W_{1}\times\cdots\times
W_{M}$, where $W_{j}$ is a space of functions on $\Omega_{j}$. The natural
splitting of the bilinear form $a\left(  \cdot,\cdot\right)  $ into integrals
over $\Omega_{j}$ is then $a\left(  u,v\right)  =\int_{\Omega}\bigtriangledown
u\cdot\bigtriangledown v=\sum_{j=1}^{N}a_{j}\left(  u,v\right)  $, where the
local bilinear forms%
\begin{equation}
a_{j}\left(  u,v\right)  =\int_{\Omega_{j}}\bigtriangledown u\cdot
\bigtriangledown v \label{mandel_mini_3_eq:local-form}%
\end{equation}
are used on $W_{j}$ instead of the bilinear form $a\left(  \cdot,\cdot\right)
$.

The space $V_{0}$ is the coarse space, and the rest of this paper deals with
its construction. It had been long understood and then formulated explicitly
in \cite{mandel_mini_3_Mandel-1990-ISS} that for condition numbers to be
independent of the number of subdomains, the coarse space needs to contain the
nullspace of the local bilinear forms $a_{j}\left(  \cdot,\cdot\right)  $. For
the scalar problem as in (\ref{mandel_mini_3_eq:local-form}), this means
constant functions, while for elasticity, the coarse space needs to contain
the rigid body modes of every substructure. Much of the development of the
coarse space in domain decomposition has been driven by the need for the
coarse space to satisfy this \emph{local nullspace condition} at the same time
as the \emph{bounded energy condition} $a\left(  v_{0},v_{0}\right)  \leq
C_{0}a\left(  v,v\right)  $, required as a part of
(\ref{mandel_mini_3_eq:energy-bound}).

\section{Some early domain decomposition methods}

\begin{figure}[ptb]
\begin{center}%
\begin{tabular}
[c]{cc}%
\includegraphics
[height=4.5cm]{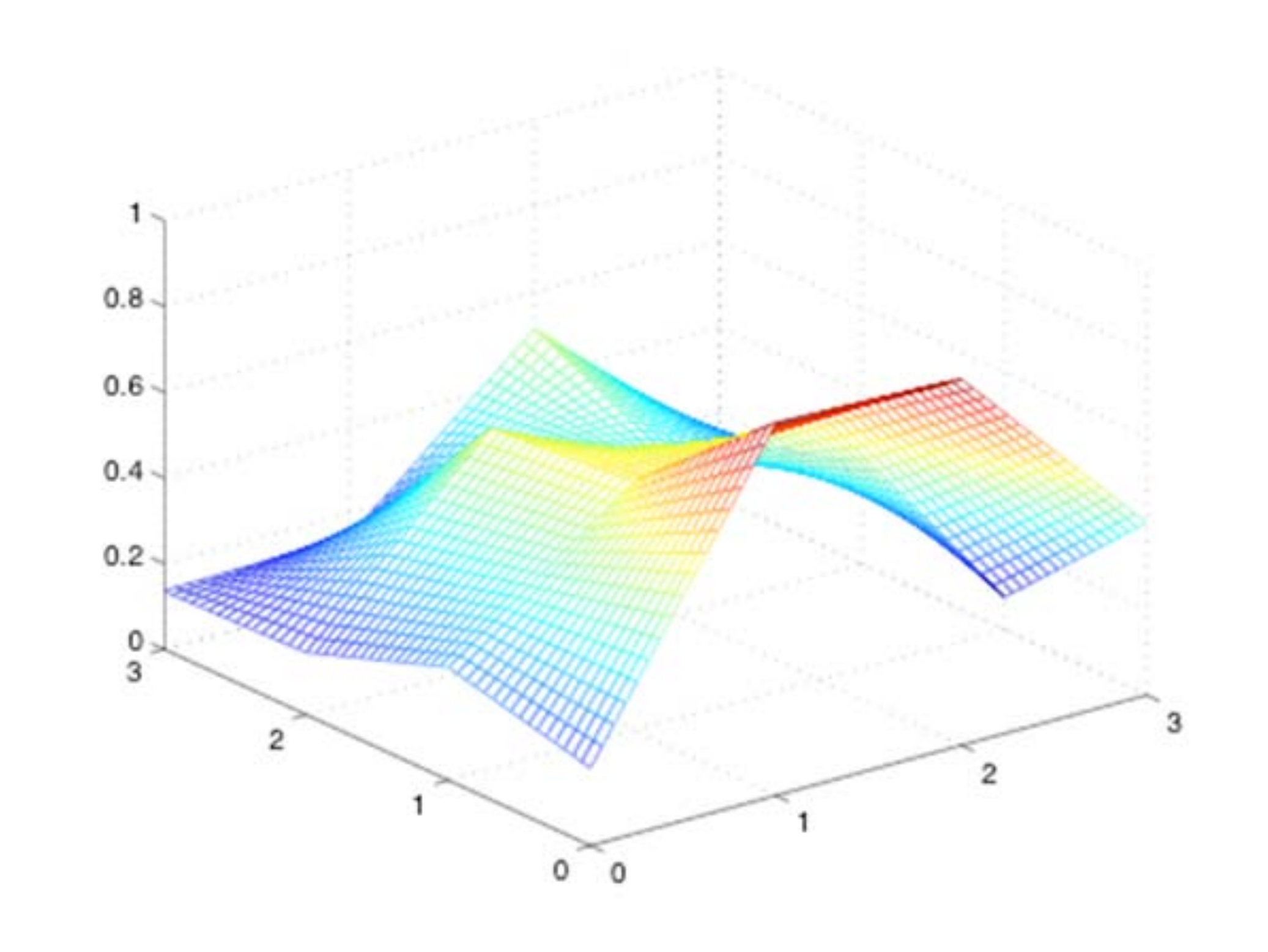} &
\includegraphics
[height=4.5cm]{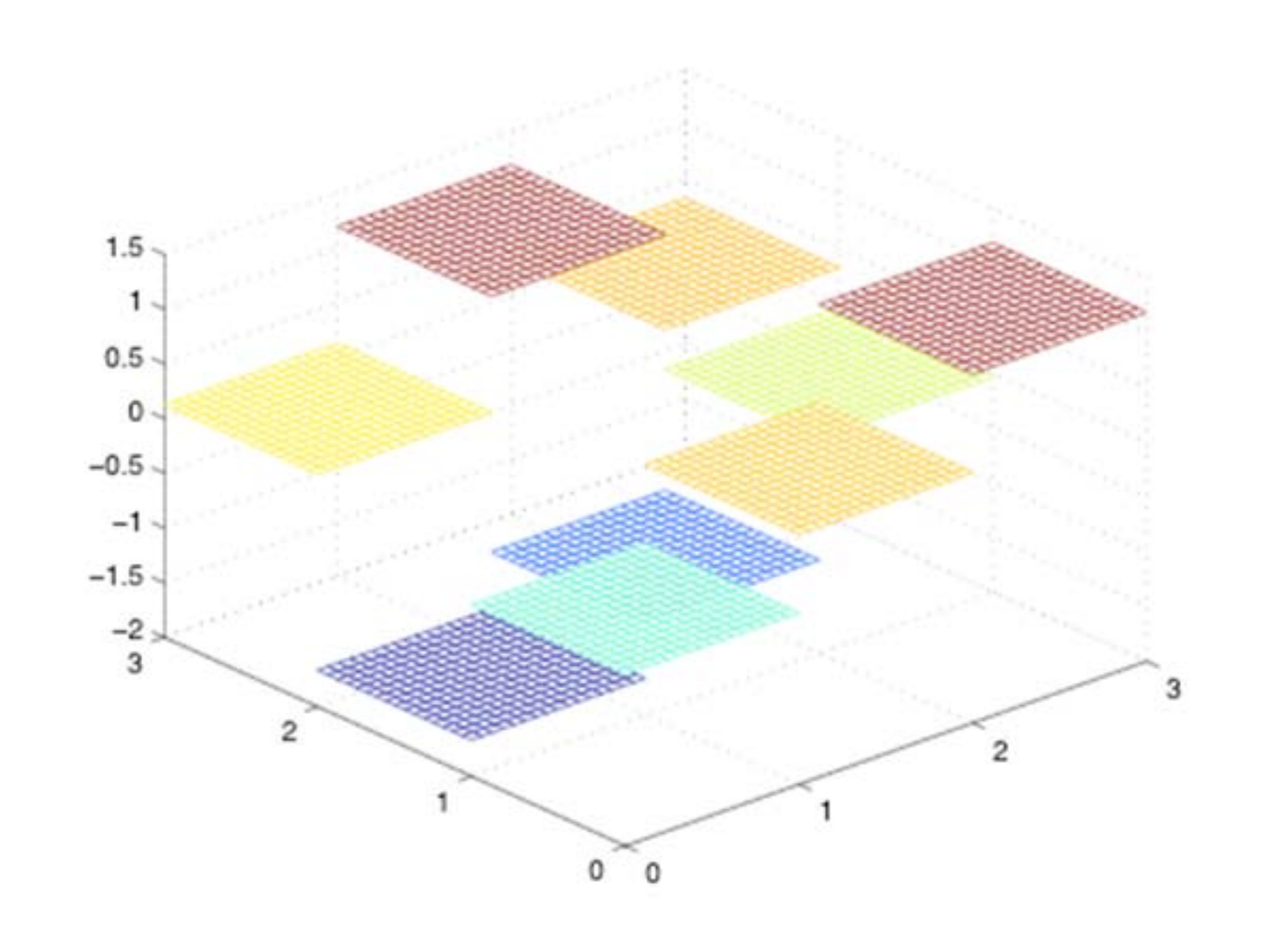}
\end{tabular}
\end{center}
\caption{Left: Piecewise bilinear coarse space function. Right: piecewise
constant functions.}%
\label{mandel_mini_3_fig:bilinear-constant}%
\end{figure}

By taking $v_{0}$ in (\ref{mandel_mini_3_eq:energy-bound}) first, we see that
the design objective of the coarse space is that there should exist a mapping
$v\in V\longmapsto v_{0}\in V_{0}$ such that (i) the energy of $v_{0}$ is not
too large, and (ii) the remainder $v-v_{0}$ can be decomposed in the spaces
$V_{i}$, $i=1,\ldots,N$, without increasing the energy too much. Definition of
$v_{0}$ by linear or bilinear interpolation is the natural first choice 
(Fig.~\ref{mandel_mini_3_fig:bilinear-constant} left). Because of the discrete
Sobolev inequality, this works fine in 2D: values of $v$ at interpolation
nodes are bounded by the energy of $v$ up to a logarithmic factor in the mesh
size $h$. The remainder $v-v_{0}$ is tied to zero by its zero values at the
interpolation nodes, and it turns out it can be decomposed into $v_{i}$'s with
bounded energy (up to a logarithmic factor). In 3D, however, the pointwise
values of $v$ for constant energy of $v$ can grow quickly as $h\rightarrow0$,
so interpolation can no longer be used. Overlapping methods
(\cite{mandel_mini_3_Dryja-1994-DDA}) use decomposition into $v_{i}$'s by a
partition of unity on overlapping subdomains $\Omega_{i}$, and they carry over
to 3D; only the interpolation from the values of $v$ needs to be replaced by a
method that is energy stable in 3D, such as interpolation from averages or
$L^{2}$ projection. In some non-overlapping methods, however, the functions
$v_{i}$ are defined in such way that they are zero at the nodes that define
the values of $v_{0}$, e.g., \cite{mandel_mini_3_Bramble-1986-CPE}. Then a
straightforward extension of the method to 3D forces $v_{0}$ to be linear
interpolant from pointwise values of $v$.
\cite{mandel_mini_3_Bramble-1989-CPE} resolved this problem by redefining the
coarse bilinear form $a_{0}$ so that $a_{0}\left(  u,u\right)  =\sum_{i=1}%
^{N}\min_{c_{i}}\int_{\Omega}\left\vert \bigtriangledown u-c_{i}\right\vert
^{2}$; cf., \cite{mandel_mini_3_Mandel-1990-TDD} for a generalization to
elasticity and an algebraic explanation. The coarse space degrees of freedom
are one number $c_{i}$ per substructure, thus the coarse space can be thought
of as piecewise constant (Fig.~\ref{mandel_mini_3_fig:bilinear-constant}
right). Piecewise constant coarse space used with the original bilinear form
$a\left(  \cdot,\cdot\right)  $ results in aggregation methods
(\cite{mandel_mini_3_Vanek-1996-AMS}). \cite{mandel_mini_3_Dryja-1988-MDD}
defined the interpolant by discrete harmonic functions, which have lower
energy than piecewise linear functions.

\section{Balancing domain decomposition (BDD) and FETI}

\begin{figure}[ptb]
\begin{center}%
\begin{tabular}
[c]{cc}%
\includegraphics[width=6cm]{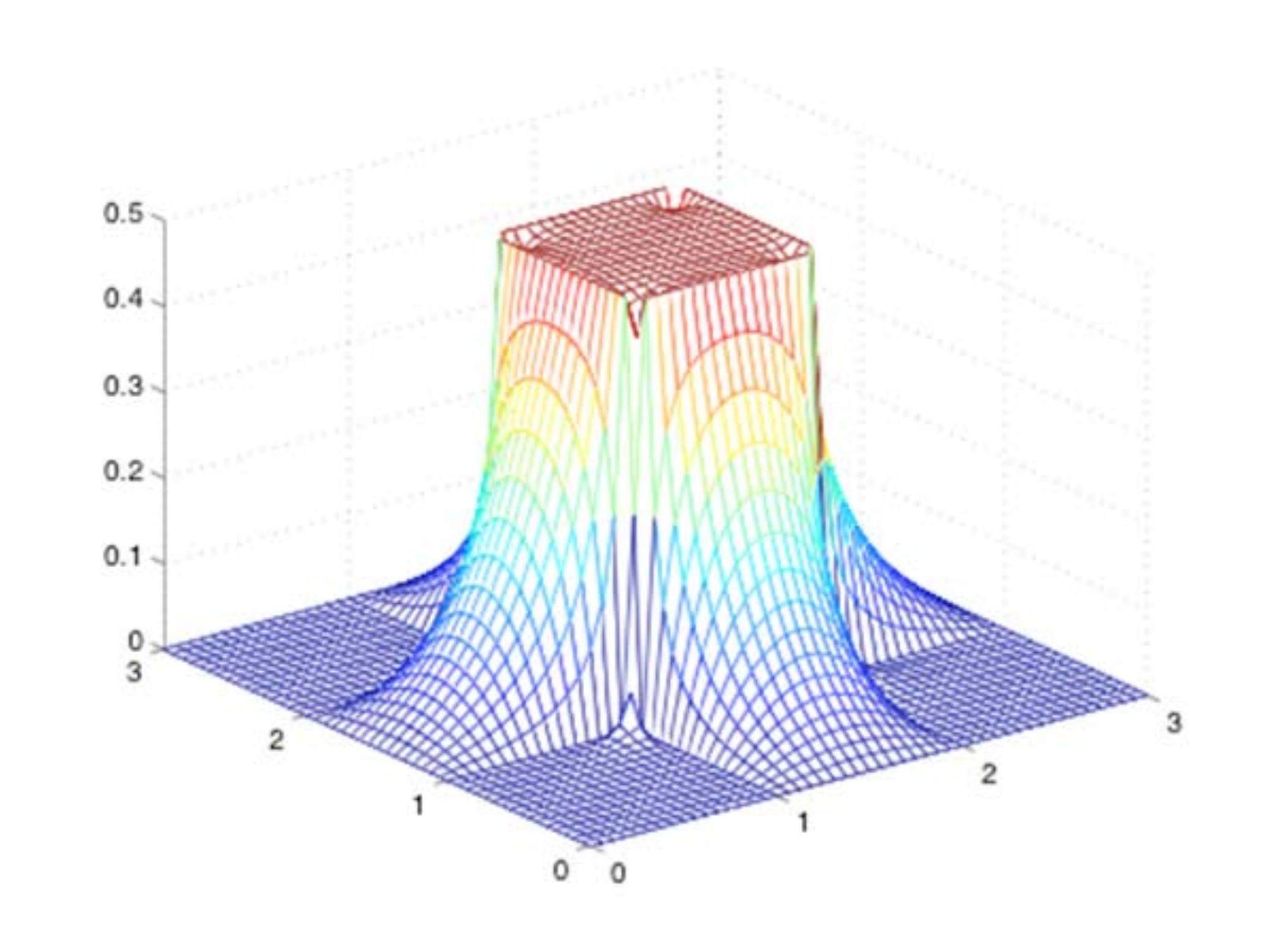} &
\includegraphics
[width=6cm]{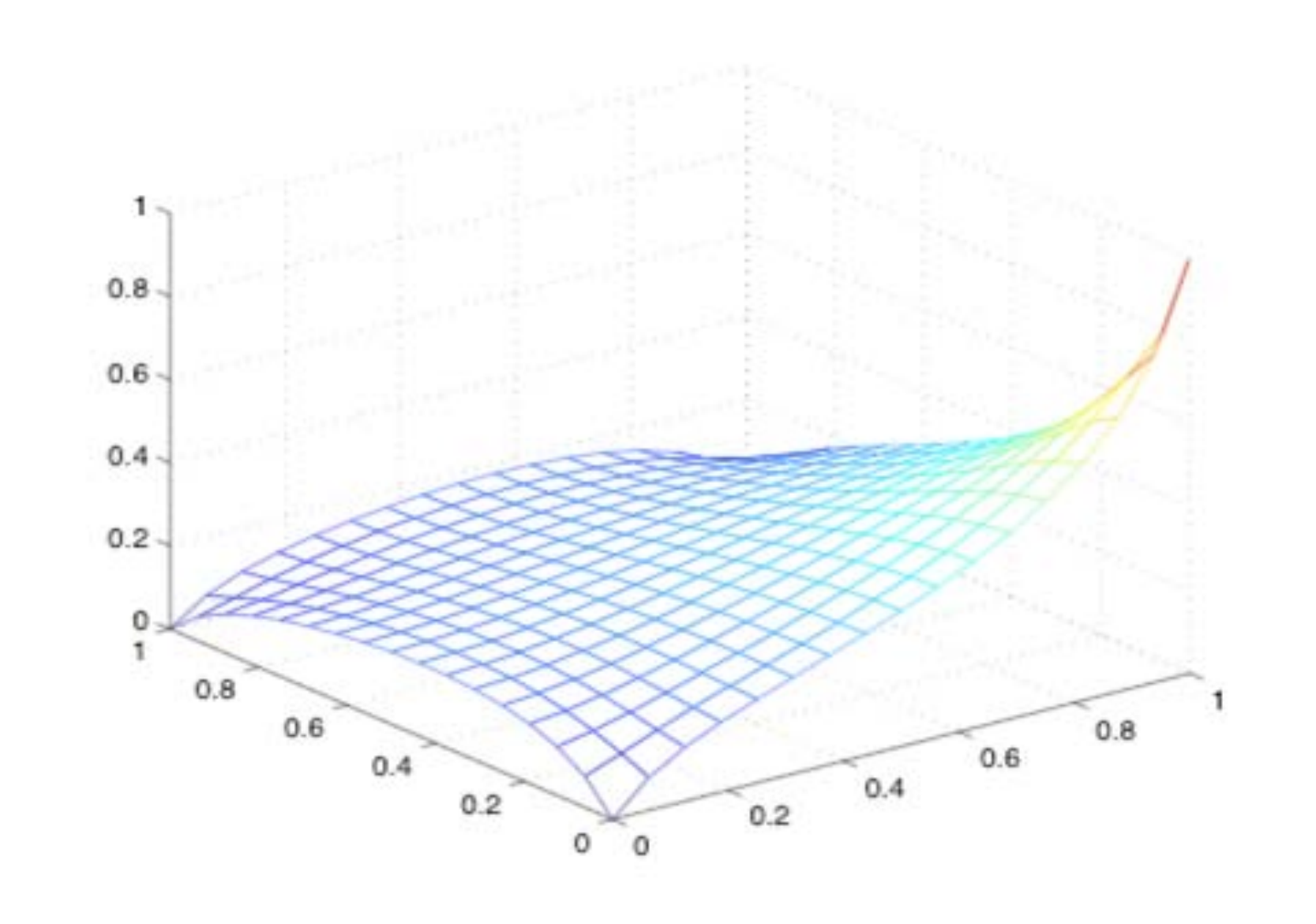}
\end{tabular}
\end{center}
\caption{Left: BDD coarse basis function, with support on one substructure and
adjacent ones. Right: Coarse function on one substructure of BDD for plates,
and BDDC (reproduced from \cite{mandel_mini_3_Mandel-2003-CBD}).}%
\label{mandel_mini_3_fig:bdd}%
\end{figure}

The BDD method was created by \cite{mandel_mini_3_Mandel-1993-BDD} by adding a
special coarse space to the Neumann-Neumann (NN) method from
\cite{mandel_mini_3_DeRoeck-1991-ATL}. The NN method uses the additive
preconditioner with the local forms $a_{i}$ from
(\ref{mandel_mini_3_eq:local-form}) and no coarse space. In the NN\ method,
the local forms are generally singular and the local problems
(\ref{mandel_mini_3_eq:subproblem}) are not consistent. The
BDD\ preconditioner applies multiplicatively a coarse correction based on a
known superspace $Z_{i}$ of the local nullspace and designed so that the
right-hand side in (\ref{mandel_mini_3_eq:subproblem}) is orthogonal to
$Z_{i}$. Since the nullspace of $a_{i}$ is contained in $Z_{i}$,
(\ref{mandel_mini_3_eq:subproblem}) is now guaranteed to be consistent. The
coarse space is obtained by averaging between adjacent substructures and
extending the functions from the substructure boundaries in the interior with
minimal energy (i.e., as discrete harmonic). A basis function of the resulting
coarse space is in Fig.~\ref{mandel_mini_3_fig:bdd} right.
Of course, for elasticity, rigid body modes are used rather than constants,
giving 6 coarse degrees of freedom per substructure in 3D.

BDD is completely \emph{algebraic}. It can be implemented only by calls to
subdomain matrix-vector multiplication and by access to a basis of the local
space $Z_{i}$ (such as the rigid body modes written in terms of the degrees of
freedom). This made possible a black-box type application of BDD to mixed
finite elements in \cite{mandel_mini_3_Cowsar-1995-BDD}: the substructure
matrix-vector multiply becomes the mapping of pressure on substructure faces
to the velocity in the normal direction.
(Some components of other methods can be generated algebraically also; e.g.,
overlapping Schwarz methods are used as smoothers in adaptive algebraic
multigrid in \cite{mandel_mini_3_Poole-2003-AAC}.)

BDD with the spaces $Z_{i}$ given by constants or rigid body modes is not suitable
for 4th order problems (such as plate bending), because the tearing at corners
has high energy - the trace norm associated with 4th order problem is the
Sobolev norm $H^{3/2}$. But empowering BDD by enriching the coarse space was
envisioned already in \cite{mandel_mini_3_Mandel-1993-BDD}, and all that was
needed was to enlarge the spaces $Z_{i}$ so that after the coarse correction,
the error is zero at corners, thus the tears across the corners do not matter.
In \cite{mandel_mini_3_LeTallec-1998-NND}, such $Z_{i}$ consists of functions
determined by their values at the corners of the substructure, and by having
minimal energy (Fig.~\ref{mandel_mini_3_fig:bdd} right).

The FETI\ method by \cite{mandel_mini_3_Farhat-1991-MFE} runs in the dual
space of Lagrange multipliers and it uses a coarse space constructed from the
exact nullspace of the local problems (\ref{mandel_mini_3_eq:subproblem}). In
the scalar case, this is the space of discontinuous piecewise constant
functions (Fig.~\ref{mandel_mini_3_fig:bilinear-constant} right), and of
piecewise rigid body modes for elasticity. Since the dual space (after
elimination of the interior) is equipped with the $H^{-1/2}$ norm, jumps
between subdomains do not cause a large energy increase. Like BDD, FETI is
completely algebraic, which is why the two methods have become popular in
practice. \cite{mandel_mini_3_Mandel-1999-SSM} generalized FETI\ to deal with
4th order problems analogously as in BDD, but the resulting method, called
FETI-2, was quite complicated. Since the basic algebra of FETI\ relies on the
exact nullspace of the local problems, the added coarse functions had to be in
a new coarse space of their own, with the additional components of the coarse
correction wrapped around the original FETI\ method. Eventually, FETI-2 was
superseded by FETI-DP.

A Neumann-Neumann method, also called balancing but somewhat different from
BDD, was developed in \cite{mandel_mini_3_Dryja-1995-SMN}. This method uses
the same coarse space as BDD, but additively, and it takes care of the
singularity in the local problems by adding small numbers to the diagonal. To
guarantee optimal condition bounds, a modification of the form $a_{0}$ is
needed.
This method is not algebraic in the same sense as BDD or FETI, i.e., relying
on the matrices only.

\section{BDDC and FETI-DP}

\begin{figure}[ptb]
\begin{center}
\hspace*{-0.6cm}
\begin{tabular}
[c]{cc}%
\includegraphics
[width=6cm]{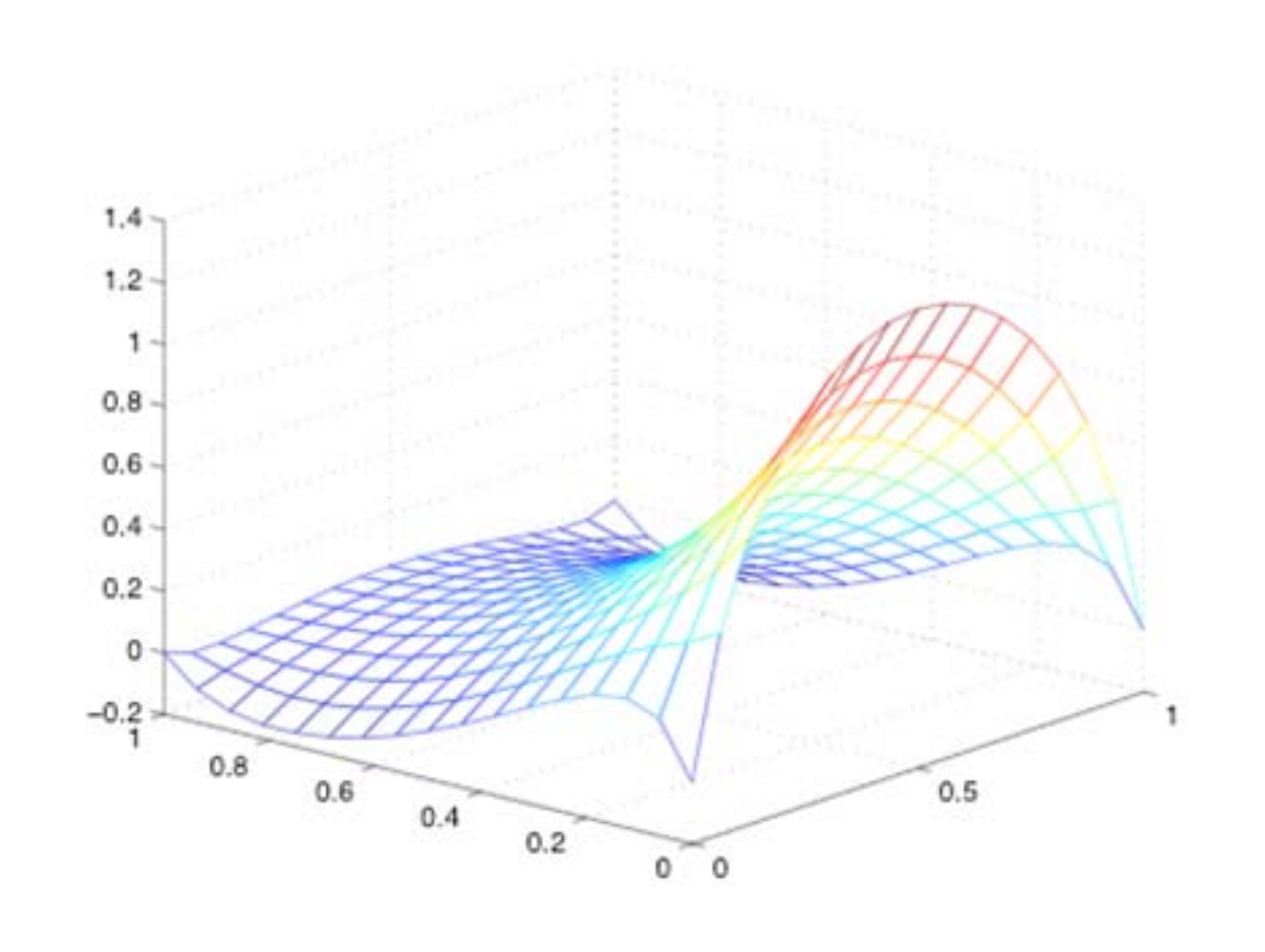} &
\includegraphics
[width=6cm]{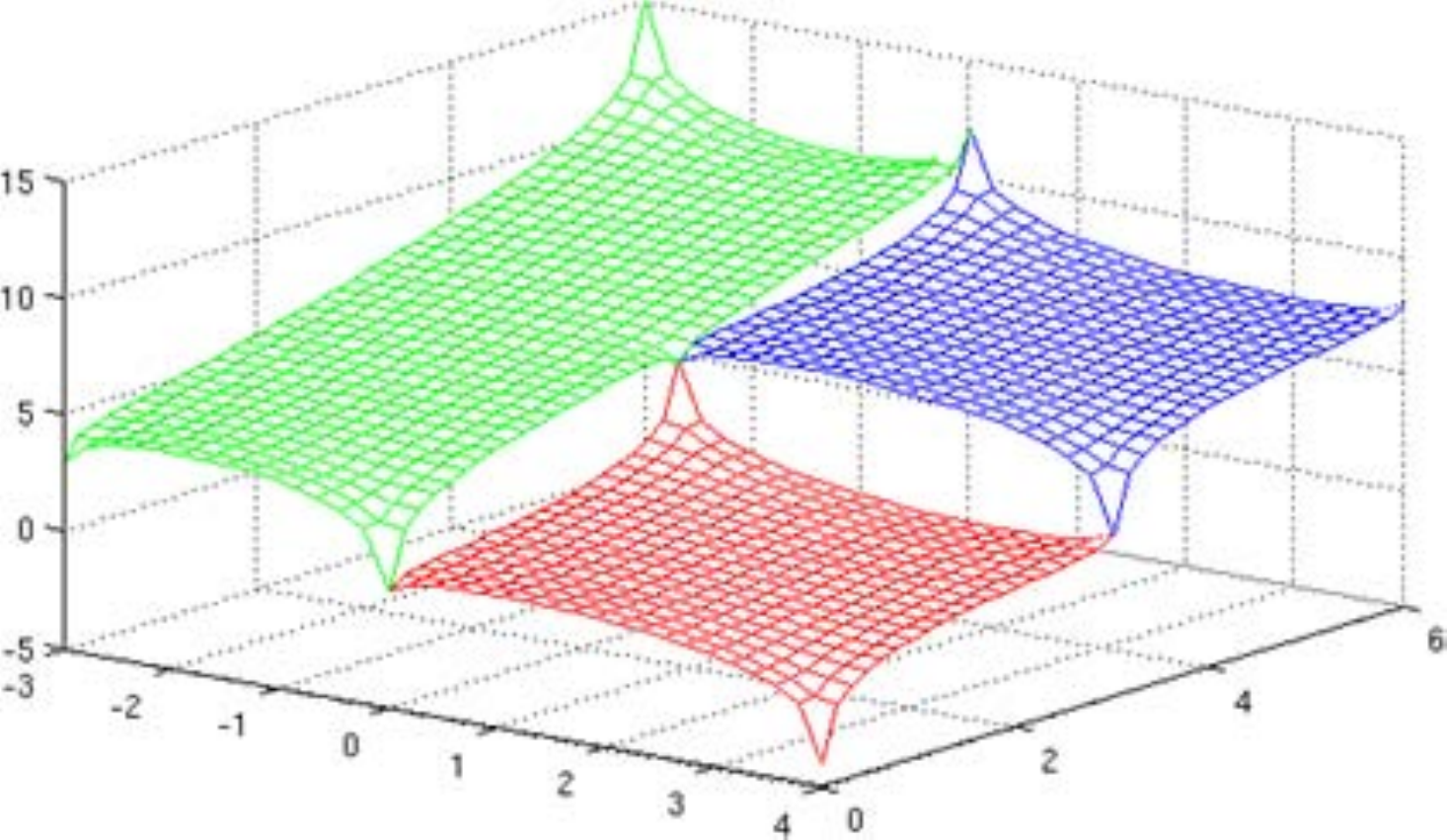}
\end{tabular}
\end{center}
\caption{Left: coarse function in BDDC for edge average degree of freedom on
one substructure. Right: BDD for plates and BDDC coarse space with vertex
degrees of freedom on several substructures (courtesy of Marta
\v{C}ert\'{\i}kov\'{a} and Jakub \v{S}\'{\i}stek).}%
\label{mandel_mini_3_fig:bddc}%
\end{figure}

A satisfactory extension of FETI and BDD to 4th order problems came only with
FETI-DP by \cite{mandel_mini_3_Farhat-2000-SDP} and BDDC by
\cite{mandel_mini_3_Dohrmann-2003-PSC}. These methods are based on identical
components and have the same spectrum. except possibly for the eigenvalues
equal to zero and one (\cite{mandel_mini_3_Mandel-2005-ATP}), so we can
discuss BDDC only. The coarse space consists of functions given by their
values of coarse degrees of freedom and energy minimal on every substructure
independently. For coarse degrees of freedom given by values on substructure
corners, this is the same coarse space as in BDD for plates in
\cite{mandel_mini_3_LeTallec-1998-NND} (Fig.~\ref{mandel_mini_3_fig:bdd}
right, Fig.~\ref{mandel_mini_3_fig:bddc} right), and the substructure spaces
$W_{i}$ are also the same. The new feature of BDDC is that the coarse
correction is additive, not multiplicative, resulting in a sparser coarse
matrix (\cite{mandel_mini_3_Mandel-2003-CBD}). In 3D, FETI-DP and BDDC require
additional degrees of freedom for optimal convergence, namely averages on
faces or edges
(\cite{mandel_mini_3_Farhat-2000-SDP,mandel_mini_3_Klawonn-2002-DPF}), cf.,
Fig.~\ref{mandel_mini_3_fig:bddc} left for a visualization in 2D.

\section{Adaptive methods by enriching the coarse space}

\label{mandel_mini_3_sec:adaptive}
\begin{figure}[ptb]
\begin{center}
\includegraphics
[width=8cm]{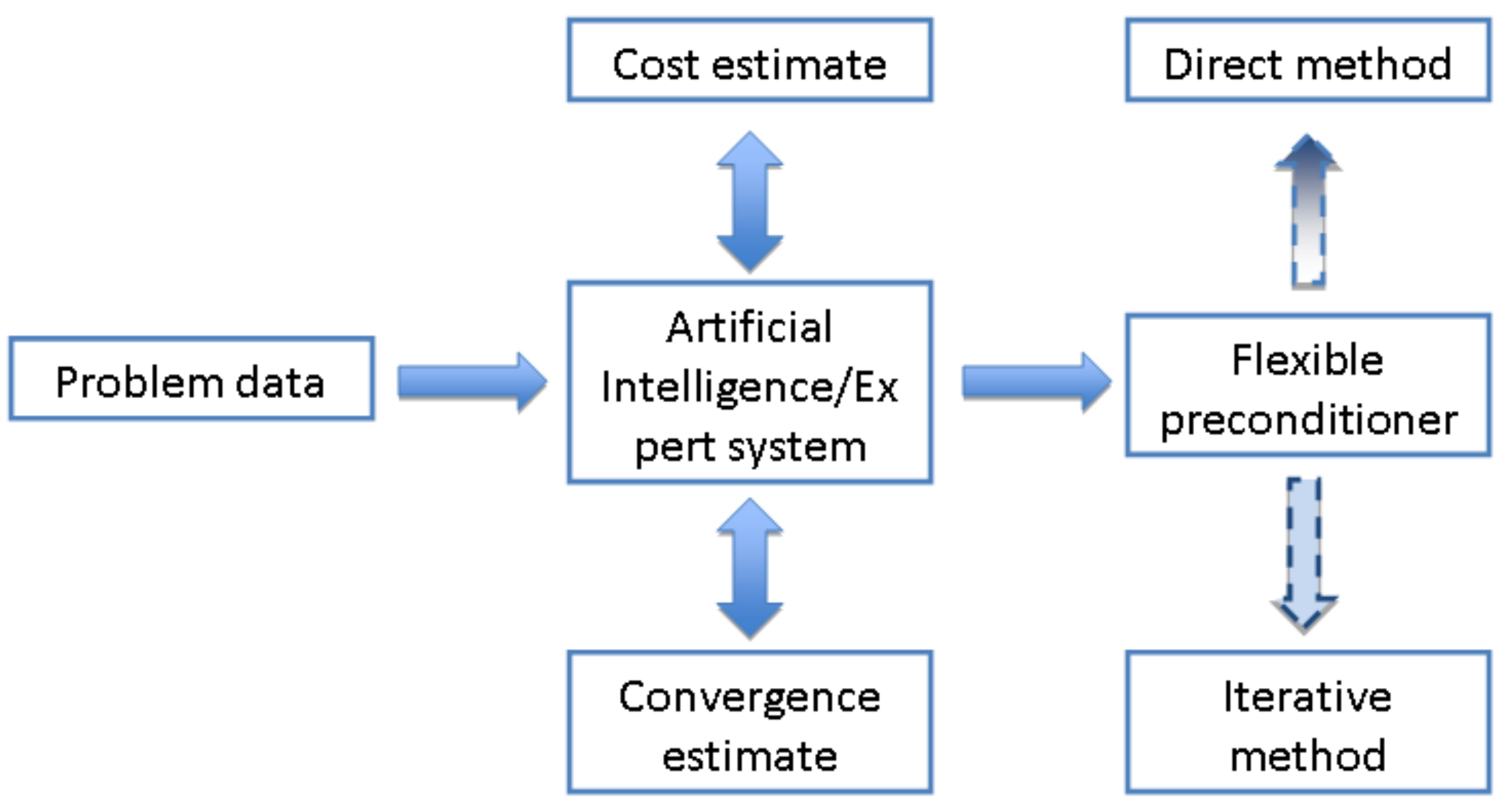}
\end{center}
\caption{Intelligent iterative method. Adapted from
\cite{mandel_mini_3_Mandel-1993-IBI}.}%
\label{mandel_mini_3_fig:intelligent}%
\end{figure}

Enlarging the coarse space is a powerful but
expensive tool. When the coarse space is the whole space, domain decomposition
turns into a direct solver. So, adding suitable functions to the coarse space
adaptively can yield a robust method, which is fast on easy problems, but does
not fail on hard ones (Fig.~\ref{mandel_mini_3_fig:intelligent}). In
\cite{mandel_mini_3_Mandel-1996-IMP}, the coarse space in the $p$-version
finite element method consists of linear functions when all is good, quadratic
functions when things get worse, all function in one direction in the case of
strong anisotropy, up to all functions when the heuristic gives up. In
\cite{mandel_mini_3_Poole-2003-AAC}, a similar methodology was applied in
algebraic multigrid. In \cite{mandel_mini_3_Mandel-2007-ASF} and in the
companion paper \cite{mandel_mini_3_Sousedik-2010-AMB} in this volume, the
coarse space in BDDC is enriched by adaptively selected linear combinations of
basis functions on substructure faces.

\bibliographystyle{plainnat}
\bibliography{mandel_mini_3}

\end{document}